\documentclass[10pt]{article}
 
\usepackage{graphicx}
\usepackage{amssymb}
\usepackage{epstopdf}
\usepackage{amsmath}
\usepackage{amsfonts}

\DeclareSymbolFont{AMSb}{U}{msb}{m}{n}
\DeclareMathSymbol{\N}{\mathbin}{AMSb}{"4E}
\DeclareMathSymbol{\Z}{\mathbin}{AMSb}{"5A}
\DeclareMathSymbol{\R}{\mathbin}{AMSb}{"52}
\DeclareMathSymbol{\Q}{\mathbin}{AMSb}{"51}
\DeclareMathSymbol{\I}{\mathbin}{AMSb}{"49}
\DeclareMathSymbol{\C}{\mathbin}{AMSb}{"43}

\begin{document}
 
\addtolength{\textheight}{0 cm}
\addtolength{\hoffset}{0 cm}
\addtolength{\textwidth}{0 cm}
\addtolength{\voffset}{0 cm}

\setcounter{secnumdepth}{5}
 \newtheorem{proposition}{Proposition}[section]
\newtheorem{theorem}{Theorem}[section]
\newtheorem{lemma}[theorem]{Lemma}
\newtheorem{coro}[theorem]{Corollary}
\newtheorem{remark}[theorem]{Remark}
\newtheorem{ex}[theorem]{Example}
\newtheorem{claim}[theorem]{Claim}
\newtheorem{conj}[theorem]{Conjecture}
\newtheorem{definition}[theorem]{Definition}
\newtheorem{application}{Application}
 
\newtheorem{corollary}[theorem]{Corollary}
\def\HADX{{\cal H}_{\rm AD}(X)}
\def\HADY{{\cal H}_{\rm AD}(Y)}
\def\HADH{{\cal H}_{\rm AD}(H)}
\def\HTADX{{\cal H}_{\rm TAD}(X)}
\def\HTADY{{\cal H}_{\rm TAD}(Y)}
\def\HTADH{{\cal H}_{\rm TAD}(H)}
\def\LX{{\cal L}(X)}
\def\LY{{\cal L}(Y)}
\def\LH{{\cal L}(H)}
 \def\ASD{{\cal L}_{\rm AD}(X)}
 \def\ASDY{{\cal L}_{\rm AD}(Y)}
\def\ASDH{{\cal L}_{\rm AD}(H)}
 \def\ASDP{{\cal L}^{+}_{\rm AD}(X)}
  \def\ASDYP{{\cal L}^{+}_{\rm AD}(Y)}
   \def\ASDHP{{\cal L}^{+}_{\rm AD}(H)}
    \def\TADX{{\cal L}_{\rm TAD}(X)}
        \def\TADY{{\cal L}_{\rm TAD}(Y)}
            \def\TADH{{\cal L}_{\rm TAD}(H)}
 \def\CX{{\cal C}(X)}
\def\CY{{\cal C}(Y)}
\def\CH{{\cal C}(H)}
 
\def\PX{{\cal A}(X)}
\def\PY{{\cal A}(Y)}
\def\PH{{\cal A}(H)}
\def\phi{{\varphi}}
\def\AH{A^{2}_{H}}
\def\B {{\cal B}} 
\def\C {{\cal C}} 
\def\H{{\cal H}}

\newcommand{\al}{\alpha}
\newcommand{\de}{\delta}

\newcommand{\ra}{\rightarrow}
\def\phi{{\varphi}}

\title{Maximal monotone operators are selfdual vector fields and vice-versa}
\author{ Nassif  Ghoussoub
\\
\small Department of Mathematics,
\small University of British Columbia, \\
\small Vancouver BC Canada V6T 1Z2 \\
\small {\tt nassif@math.ubc.ca} 
\\
}
\maketitle

\begin{abstract} 
\noindent If $L$ is a selfdual Lagrangian $L$ on a reflexive phase space  $X\times X^*$,  then the vector field $x\to \bar\partial L(x):=\{p\in X^*; \, (p,x)\in \partial L(x,p)\}$ is maximal monotone. Conversely, any maximal monotone operator $T$ on  $X$ is derived from such a potential on phase space, that is there exists a selfdual Lagrangian $L$ on $X\times X^*$ (i.e, $L^*( p, x) =L(x, p )$) such that $T=\bar\partial L$. 

 This solution to  problems raised by Fitzpatrick can be seen as an extension of a celebrated result  of Rockafellar stating that  maximal cyclically monotone operators are actually of the form $T=\partial \phi$ for some convex lower semi-continuous function on $X$.  This representation  allows  for the application of the selfdual variational theory --recently developed by the author-- to the equations  driven by maximal monotone vector fields. Consequently, solutions to  equations of the form $\Lambda x\in Tx$ for a given map $\Lambda: D(\Lambda)\subset X\to X^*$, can now be obtained by minimizing functionals of the form $I(x)=L(x,\Lambda x)-\langle x, \Lambda x\rangle$.

\end{abstract}

\section{Introduction} In a series of papers (\cite{G2}--\cite{GT2}), we developed and used   the concept of selfdual Lagrangians $L$ on phase space $X\times X^*$ --were $X$ is a Banach space-- in order  to provide variational formulations and resolutions to various differential equations which are not variational in the sense of Euler-Lagrange. The main idea behind this theory --which is developed in its full generality in the upcoming monograph \cite{G10}--  originated from the fact that a large set of PDEs and evolution equations can be written in the form 
\begin{equation} \label{basic.equation}
(p,x)\in \partial L(x,p),
\end{equation}
where $\partial L$ is the subdifferential of a   Lagrangian $L:X\times X^*\to \R\cup \{+\infty\}$ that is convex  and lower semi-continuous --in both variables-- while satisfying 
\begin{equation}
\hbox{$L^*( p, x) =L(x, p )$   for all $(p,x)\in X^{*}\times X$.}
\end{equation}
Here $L^*$ is the Legendre transform of $L$ in both variables, that is
\[
L^*(p,x)=\sup\left\{ \langle p,y\rangle +\langle x,q\rangle -L(y, q);\, (y, q)\in X\times X^*\right\}.             
\]
Solutions are then found for a given $p$ by simply minimizing the functional $I_p(x)=L(x,p)-\langle x,p\rangle$ and by proving that the minimum is actually zero. In other words, by defining the {\it derived vector fields of L} at $x\in X$ 
 to be the --possibly empty-- sets
\begin{equation} \label{selfdual.vector}
\hbox{${\overline\partial} L (x):=\{p\in X^*; \, L(x,p)-\langle x,p\rangle =0\}=\{p\in X^*; \, (p,x)\in \partial L(x,p)\}$,}
\end{equation}
one can then find variationally the zeroes of  those set-valued maps $T: X\to 2^{X^*}$ of the form $T(x)=\bar\partial L(x)$ for some selfdual Lagrangian $L$ on $X\times X^*$. 

These {\it selfdual vector fields} are 
natural extensions of  subdifferentials of convex lower semi-continuous functions. Indeed, 
 the most basic selfdual Lagrangians are of the form $L(x,p)=\phi (x) +\phi^*(p)$ where $\phi$ is  such a  function on $X$, and $\varphi^{*}$ is its Legendre conjugate on $X^{*}$,  in which case ${\overline \partial}  L (x)=\partial \phi (x)$.
 The corresponding variational problem (i.e., minimizing $I(x)=L(x,0)=\phi (x) +\phi^*(0)$) reduces to the classical approach of minimizing a convex functional in order to solve equations of the form $0\in \partial \phi (x)$. 

More interesting examples of selfdual Lagrangians are of the form 
$
L(x,p)=\varphi (x) +\varphi^{*}(-\Gamma x+p)
$
 where $\varphi$ is a convex and lower semi-continuous function on $X$, and  $\Gamma:X\to X^{*}$ is a skew-symmetric operator. The corresponding selfdual vector field is then, 
\[
\hbox{${\overline \partial} L (x)=\Gamma x+\partial \phi (x)$.}
\]
More generally,  if the operator $\Gamma$ is merely non-negative (i.e., $\langle \Gamma x, x\rangle \geq 0$), then one can still write the vector field  $\Gamma +\partial \phi$ as ${\overline \partial} M$ for some selfdual Lagrangian $M$ defined now on $X\times X^*$, as
\[
M(x,p)=\psi (x) +\psi^{*}(-\Gamma^{as} x+p)
\]
where $\psi$ is the convex function $\psi (x)=\frac{1}{2} \langle \Gamma x, x\rangle +\varphi (x)$,  and where $\Gamma^{as}=\frac{1}{2}(\Gamma -\Gamma^{*})$ is the anti-symmetric part of $\Gamma$, and $\Gamma^{sym}=\frac{1}{2}(\Gamma +\Gamma^{*})$ is its symmetric part. The main interest being that equations of the form  $\Gamma x+\partial \phi(x)=p$ can now be solved for a given  $p\in X^*$, by simply minimizing the functional 
\[
I_p(x)=\psi (x) +\psi^{*}(\Gamma^{as} x+p)-\langle x, p\rangle
\]
 and proving that its infimum is actually zero.

More classical is the fact that the above examples are typical  set-valued nonlinear operators $T$ that  are {\it monotone}, meaning that their graphs $G(T)=\{ (x,p)\in X\times X^*; p\in T(x)\}$ are {\it monotone} subsets of 
$X\times X^*$, i.e., it satisfies:
\begin{equation}
\hbox{$\langle x-y, p - q \rangle \geq 0$ for every $(x,p)$ and $(y, p)$ in $G(T)$.}
\end{equation}
They are actually  {\it maximal monotone operators} meaning that their graph $G(T)$ is maximal in the family of monotone subsets of $X\times X^*$, ordered by set inclusion. For a comprehensive description of this theory, we refer to the monograph of Brezis \cite{Br}.

It was therefore a natural question whether one can associate a selfdual Lagrangian $L$ to any maximal operator $T$, so that equations involving such operators can be resolved variationally.  We recently realized that Krauss \cite{K} and Fitzpatrick \cite{F}  had done some work in this direction in the 80's, and have managed to associate to a maximal monotone operator $T$, a ``sub-selfdual Lagrangian", i.e. a convex lower semi-continuous function $L$ on state space $X\times X^*$ satisfying $L^*(p,x) \geq L(x,p)$ on $ X\times X^*$, and in such a way that $T=\bar \partial L$.  The main interest of this note is to show that one can actually construct a selfdual Lagrangian associated to $T$, which answers one of the questions in \cite{F}.

 \begin{theorem}  \label{ghoussoub} Let $L$ be a proper selfdual Lagrangian $L$ on a reflexive Banach space $X\times X^*$, then the vector field $x\to \bar \partial L (x)$ is maximal monotone.

 Conversely, if $T: D(T) \subset X\to 2^{X^*}$ is a maximal monotone operator with a non-empty domain, then there exists a selfdual Lagrangian $N$ on $X\times X^*$ such that $T= \bar \partial N$. 
 
  \end{theorem}

It is  worth comparing this result with the following celebrated result of Rockafellar \cite{Ph}, which gives an integral representation of those maximal monotone operators $T: X\to 2^{X^*}$ that are {\it cyclically monotone}, i.e., those that verify for any finite number of points $(x_i, p_i)_{i=0}^n$ in the graph $G(T)$ with $x_0=x_n$, we have 
\begin{equation}
\hbox{$\sum\limits_{i=1}^n\langle p_k, x_k-x_{k-1} \rangle \geq 0$.}
\end{equation}

\begin{theorem} [Rockafellar] If $\phi:X\to \mathbb{R}\cup\{+\infty\}$ is a proper convex and lowersemi
continuous functional on a Banach space $X$, then its subdifferential map  $x\to \partial \phi (x)$ is a maximal cyclically monotone map.

 Conversely if $T:X\to 2^{X^*}$ is a maximal cyclically monotone map with a non-empty domain, then there exists a a proper convex and lower semi continuous functional on  $X$ such that $T=\partial \phi$.
 \end{theorem}

The advantages of identifying  maximal monotone operators with the class of selfdual vector fields are numerous:

\begin{enumerate}
\item As mentioned above, all equations, systems, variational inequalities, and  dissipative initial value  parabolic problems which traditionally involve maximal monotone operators,  can now be formulated and resolved variationally.   These problems can therefore be analyzed with the full range of methods --computational or not-- that are available for variational settings.  

\item  While issues around the superposition of,  and other operations on, maximal monotone operators are often delicate to prove,  the class of selfdual Lagrangians  possesses remarkable  permanence properties that are also easy to establish.   It reflects most variational aspects of convex analysis, and is stable under similar type of operations making the calculus of selfdual Lagrangians as manageable, yet much more encompassing  \cite{G2}.

\item Selfduality allows for the superposition of appropriate boundary Lagrangians  with ``interior"   Lagrangians,  leading to the resolution of problems with various linear and nonlinear boundary constraints that are not amenable to the standard variational theory \cite{GT2}. They can also be iterated with certain nonlinear operators that are far from being maximal monotone \cite{G3}.

\item Selfdual Lagrangians defined on state spaces also lift to selfdual Lagrangians on path spaces leading to a unified approach for stationary and dynamic equations.  More precisely,  selfdual  flows of the form  $\dot u (t) \in \overline \partial L(u(t))$ with a variety of time-boundary conditions can be reformulated and resolved  as $0\in \overline \partial {\cal L} (u)$ where ${\cal L}$ is a corresponding selfdual Lagrangian on path space, a phenomenon that leads to  natural and quite interesting iterations (See \cite{G2}, \cite{GT3}).  

\item The class of $B$-selfdual Lagrangians corresponding to  a general automorphism $B$ and their $B$-selfdual vector fields  goes beyond the theory of maximal monotone operators. They also cover the theory of Hamiltonian systems --when $B$ is the symplectic matrix $J$-- and other twisted  differential operators  \cite{GM2}. In return, the connection between the two concepts leads to a rich parallel theory of $B$-monotone operators which coincides with the classical theory when $B$ is the identity. The corresponding  class of PDEs and their variational principles are studied in detail in the upcoming monograph \cite{G10}.  

\end{enumerate}

 \section{Selfdual vector fields are maximal monotone operators}
 
 To any convex lower semi-continuous Lagrangian $L:X\times X^* \to \R \cup\{+\infty\}$, we can associate  a vector field $ \delta L$ at any $x\in X$,  
 to be the --possibly empty-- set
\begin{equation} \label{vector,field}
\hbox{${\delta} L (x):=\{p\in X^*; \, (p,x)\in \partial L(x,p)\}$.}
\end{equation}
It is easy to see that the convexity of $L$ yields that $x\to \delta L (x)$ is a monotone map, since if $p\in \delta L(x)$ and $q\in \delta L(y)$, then
\[
\langle p-q, x-y\rangle =\frac{1}{2}\langle (p,x)-(q,y), (x,p)-(y,q) \rangle \geq 0.
\]
We can also associate another --not always necessarily monotone-- vector field as: 
\begin{equation} \label{selfdual.vector}
\hbox{${\overline\partial} L (x):=\{p\in X^*; \, L(x,p)-\langle x,p\rangle =0\}.$}
\end{equation}

Note that $\bar \partial L$ should not be confused with the subdifferential $\partial L$ of $L$ as a convex function on $X\times X^*$. If now $L$ satisfies 
\begin{equation}\label{positive}
\hbox{$L(x,p)\geq \langle x, p\rangle$ for all $(x,p)\in X\times X^*$},
\end{equation}
then it is easy to see that $\bar \partial L (x) \subset \delta L (x)$. Moreover, if  $L$ is selfdual, i.e., if
\begin{equation} \label{subselfdual}
\hbox{$L^*(p,x)= L(x,p)$ for all $(x,p)\in X\times X^*$},
\end{equation}
then the reverse inclusion is also true since then Legendre-Fenchel duality applied to  $L$ yields:
 \[
 2L^*(p,x)=2L(x,p)= L^*(p, x)+L(x, p) \geq  2\langle x, p\rangle  
  \]
 and therefore $L(x, p)-\langle x, p\rangle =0$ if and only if $L(x, p)+L^*(p, x) =   \langle x, p\rangle +\langle p, x\rangle$ which is equivalent to $(p, x)\in \partial L(x,p)$. In other words, we can state the following observation. 

\begin{lemma} If $L$ is a convex lower semi-continuous selfdual Lagrangian on $X\times X^*$, then $p\in {\overline \partial}  L (x)$ if and only if $p\in \delta L (x)$.
\end{lemma}
Also note that  $p\in {\overline \partial}  L (x)$ if and only if  $0\in {\overline \partial}  L_p (x)$ where $L_p$ is the selfdual Lagrangian $L_p(x, q)=L(x, p+q)+\langle x, p\rangle$. This is also equivalent  to the statement that the infimum of the functional $I_p(u)=L(u, p)-\langle u, p\rangle$ is zero and is attained at $x\in X$. This leads to the following proposition which is a particular case of a more general result established in \cite{G2}.

\begin{proposition} \label{main.existence} Let $L$ be a convex lower semi-continuous  selfdual Lagrangian  on a reflexive Banach space $X\times X^{*}$, such that  for some $x_0\in X$, the function $p\to L(x_0,p)$ is  bounded on the balls of $X^*$.   Then for each $p\in X^*$, there exists $\bar x\in X$ such that: 
 \begin{equation}
 \left\{ \begin{array}{lcl}
\label{eqn:main.existence}
  L( \bar x, p)-\langle \bar x, p\rangle&=&\inf\limits_{x\in X}\left\{ L(x,p)-\langle x, p\rangle\right\}=0.\\
 \hfill  p &\in & \bar \partial L (\bar x).
\end{array}\right.
 \end{equation}
 \end{proposition}
\noindent{\bf Proof:} We can assume that $p=0$ by considering the translated  Lagrangian  $M(x, q)=L(x, p+q)-\langle x, p\rangle$ which is also selfdual on $X\times X^*$. In this case $M(x,0) \geq 0$ for every $x\in X$.

 Now consider  $({\cal P}_{q})$ the primal  minimization problem
$h(q)=\inf\limits_{x\in X}M(x,q)$
in such a way that $({\cal P}_{0})$ is the initial problem $h(0)=\inf\limits_{x\in X}M(x,0)$, 
and the dual problem $({\cal P}^{*})$ is therefore
$\sup_{y\in X}-M^{*}(0,y)$.  

We readily have the following weak duality formula:
\[
\inf {\cal P}_{0}:=\inf_{x\in X}M(x,0) \geq 0\geq \sup_{y\in X}-M^{*}(0,y):=\sup{\cal P}^{*}. 
\]
Note that  $h$ is convex on $X^{*}$, and that its Legendre conjugate satisfies  for all $y\in X$.
\[
h^{*}(y)=M^{*}(0,y)=M(y,0)=L(y, p)-\langle y, p\rangle
\]
Moreover,  $h(q)=\inf\limits_{x\in X}M(x,q) \leq L(x_0, p+q)-\langle x_0, p\rangle$ and therefore $q\to h(q)$   is bounded above on the balls of  $X^{*}$, 
and hence    
it is subdifferentiable at $0$ (i.e., the problem $({\cal P}_{0})$ is then stable). Any point $\bar x \in \partial h (0)$ satisfies $h(0) +h^{*}(\bar x)=0$, which means that
\[
-\inf_{x\in X}M(x,0)=-h(0)=h^{*}(\bar x)=M^{*}(0, \bar x)=M(\bar x,0)\geq  \inf_{x\in X}M(x,0).
\]
It follows that $\inf_{x\in X}M(x,0)=M(\bar x,0)\leq 0$ and the infimum of $({\cal P})$  is therefore zero and is attained at $\bar x$.  

\begin{lemma} \label{inf.convolution} Let $L$ be a proper selfdual Lagrangian on $X\times X^*$, then for any  convex continuous function $\phi$ on $X$, the Lagrangian defined by
\[
M(x,p)=\inf\left\{L(x, p-r)+\phi (x) +\phi^*(r); \, r\in X^*\right\}
\]
is also selfdual on $X\times X^*$.
\end{lemma}
 {\bf Proof:} Indeed,  fix $(q, y)\in X^{*}\times X$ and write:
\begin{eqnarray*}
 M^{*} (q,y)
  &=&\sup\{\langle q, x\rangle +  \langle y, p\rangle- L(x, p-r)- \phi (x) - \phi^*(r); (x, p,r)\in X\times X^*\times X^{*}\}\\
&=&  \sup\{\langle q, x\rangle +  \langle y, r+s\rangle- L(x, s)- \phi (x) - \phi^*(r); (x, s, r)\in X\times X^*\times X^{*}\}\\
&=&\sup_{x\in X} \left\{\langle x,q\rangle +\sup\limits_{(s, r)\in X^*\times X^*}\{\langle y, r+s\rangle -L(x,s)- \phi^*(r)\}- \phi (x)\right\}\\
&=&\sup_{x\in X}\left\{\langle x,q\rangle +\sup_{s\in X^*}\{\langle y,s\rangle -L(x,s)\} +\sup_{r\in X^*} \{\langle y,r\rangle -\phi^* (r)\}-\phi (x)\right\}\\
&=&\sup_{x\in X}\left\{\langle x,q\rangle +\sup_{s\in X^*}\{\langle y,s\rangle -L(x,s)\} + \phi (y)-\phi (x)\right\}\\
&=&\sup_{x\in X}\sup_{s\in X^*}\left\{\langle x,q\rangle +\langle y,s\rangle  -L(x,s)-\phi (x)\right\} + \phi (y)\\
&=& (L+T_\phi)^{*}(q,y)+\phi (y)
\end{eqnarray*}
where $T_\phi (x,s):=\phi (x) $ for all $(x,s)\in X\times X^{*}$. Note now that 
\begin{eqnarray*}
T_\phi^*(q,y)=\sup_{x,s}\left\{\langle q,x\rangle +\langle y,s\rangle -  \phi (x)\right\}=\left\{\begin{array}{lll}+\infty &\hbox{if }&y\ne 0\\ \phi^*(q)&\hbox{if }&y=0\end{array} \right.
\end{eqnarray*}
in such a way that by using the duality between sums and convolutions in both variables, we get
\begin{eqnarray*}
(L+T_\phi)^{*}(q,y)&=&{\rm conv}(L^*,T_\phi^*)(q,y)\\
&=&\inf_{r\in X^{*},z\in X}\left\{ L^*(r,z)+T_\phi^*(-r+q,-z+y)\right\}\\
&=&\inf_{r\in X^*}\left\{ L^*(r,y)+ \phi^*(-r+q)\right\}
\end{eqnarray*}
and  finally 
\begin{eqnarray*}
M^{*} (q,y)&=&(L_\phi+T)^{*}(q,y)+\phi (y)\\
&=&\inf_{r\in X^*}\left\{ L^*(r,y)+\phi^*(-r+q)\right\}+\phi (y))\\
&=&\inf_{s\in X^*}\left\{ L^*(q-s,y)+\phi^*(s)\right\}+\phi (y))\\
&=&\inf_{s\in X^*}\left\{ L(y,q-s)+\phi (y) +\phi^*(s) \right\}\\
  &=& M(q,y).
\end{eqnarray*}

 \begin{proposition} \label{subdifferential} Let $L$ be a selfdual Lagrangian $L$ on a reflexive Banach space $X\times X^*$. The following assertions then hold:
  \begin{enumerate}
  \item The vector field $x\to \bar \partial L (x)$ is maximal monotone.
 
  \item  If $L$ is strictly convex in the second variable, then the maximal monotone vector field  $x\to \bar \partial L (x)$ is single-valued on its domain.
   
   \item If $L$ is uniformly convex in the second variable (i.e.,  if  $L(x,p)-\epsilon\frac{\| p\|^2}{2}$  is convex in $p$ for some $\epsilon>0$) 
    then the vector field  $x\to \bar \partial L (x)$ is a Lipschitz  maximal monotone operator on its domain. 
  \end{enumerate}
    \end{proposition}   
{\bf Proof:} Denoting by $J:X\to 2^{X^*}$  the duality map between $X$ and $X^*$, it suffices to show that the vector field $\bar\partial  L+    J$ is onto \cite{Ph}.   In other words, we need to find for any $p\in X^*$,  an $x\in  X$ such that $p=\bar\partial  L(x)+  J(x)$.
For that, we consider the following Lagrangian on $X\times X^*$.
\[
M(x,p)=\inf\left\{L(x, p-r)+\frac{1}{2}\|x\|_X^2 + \frac{1}{2}\|r\|_{X^*}^2; \, r\in X^*\right\}.
\]
It is a selfdual Lagrangian according to the previous lemma. Moreover -- assuming without loss of generality-- that the point $(0,0)$ is in the domain of $L$, we get   the estimate 
$
  M(0,p)\leq L(0,0)+\frac{1}{2}\|p\|_{X^*}^2  
$, and therefore Proposition \ref{main.existence} applies and we obtain $\bar x \in X$ so that $ p\in \bar\partial  M ( \bar x)$.
This means that 
\[
M(\bar x,p)-\langle \bar x, p\rangle=\inf\left\{L(\bar x, p-r)-\langle \bar x, p-r\rangle+\frac{1}{2}\|\bar x\|_X^2 + \frac{1}{2}\|r\|_{X^*}^2-\langle \bar x, r\rangle;\,  r\in X^* \right\}=0.
\]
In other words, there exists $\bar r\in  J(\bar x)$ such that  $p - \bar r  \in  \bar\partial  L ( \bar x)$ and we are done.  

The other assertions of the proposition are straightforward and left to the interested reader. 
 
\section{Maximal monotone operators are selfdual vector fields}

We start with the following lemma which is essentially due to Fitzpatrick \cite{F}. 

\begin{lemma} \label{fitz}Let $T: D(T) \subset X\to 2^{X^*}$ be a monotone operator, and consider  on $X\times X^*$ the Lagrangian $L_T$ defined by
\begin{equation}
L_T(x,p)=\sup\{\langle p,y\rangle +\langle q, x-y\rangle; \, (y,p)\in G(T)\}
\end{equation}
\begin{enumerate}
\item If ${\rm Dom} (T)\neq \emptyset$, then $L_T$ is convex and lower semi-continuous function on $X\times X^*$ such that for every $x\in {\rm Dom} (T)$, we have $Tx\subset \bar \partial L(x) \cap \delta L (x)$. Moreover, we have 
\begin{equation}
\hbox{$L_T^*(p,x) \geq L_T(x,p)$ for every $(x,p)\in X\times X^*$.}
\end{equation}

\item If $T$ is maximal monotone, then $T= \bar \partial L=\delta L$
and 
\begin{equation} \label{fenchel}
\hbox{$L_T(x,p)\geq \langle x,p\rangle$ for all $(x,p)\in X\times X^*$.} 
\end{equation}
\end{enumerate}
\end{lemma} 

{\bf Proof:} (1) If $x\in {\rm Dom} (T)$ and  $p\in Tx$, then the monotonicity of $T$ yields for any $(y, q)\in G(T)$
\[
\langle x,p\rangle \geq \langle y,p\rangle +\langle x-y, q\rangle
\]
in such a way that $L_T(x,p)\leq \langle x,p\rangle$. On the other hand, we have 
\[
L_T(x,p) \geq \langle x,p\rangle +\langle p, x-x\rangle =\langle x,p\rangle,
\]
and therefore $p\in \bar \partial L(x)$.

Write now for any $(y, q)\in X\times X^*$, 
\begin{eqnarray*}
L_T(x+y, p+q)-L_T(x,p)&=&\sup\left\{ \langle p+q, z\rangle +\langle r, x+y\rangle-\langle z,r\rangle;\, (z,r)\in G(T)\right\}-L_T(x,p)\\
&\geq &\langle p+q, x\rangle +\langle p, x+y\rangle-\langle p,x\rangle -\langle p,x\rangle\\
&= &\langle q, x\rangle +\langle p, y\rangle 
\end{eqnarray*}
which means that $(p,x) \in \partial L (x,p)$ and therefore $p\in \delta L(x)$.

Note also that $L_T(x,p)=H_T^*(p,x)$ where $H_T$ is the Lagrangian on $X\times X^*$ defined by 
\begin{equation} \label{hamiltonian}
 H_T(x,p)= \left\{ \begin{array}{lcl}
\hbox{$\langle x, p \rangle$ \quad if $(x,p)\in G(T)$}\\
\hbox{$+\infty$ \quad \quad \quad otherwise.} 
\end{array}\right.
\end{equation}
Since $L_T(x, p)=\langle x, p \rangle=H_T(x,p)$ whenever $(x,p)\in G(T)$, it follows that 
$L_T\leq H_T$ on $X\times X^*$ and so $L_T^*(p,x)\geq H_T^*(p,x)=L_T(x,p)$ everywhere. \\

(2) If now $T$ is maximal then necessarily $Tx= \delta L(x)\cap \bar \partial L(x)=\delta L(x)$ since $x\to \delta L(x)$ is a monotone extension of $T$.  

In order to show (\ref{fenchel}), assume to the contrary that $L_T(x,p)< \langle x,p\rangle$ for some $(x,p)\in X\times X^*$. It follows that 
\[
\hbox{$\langle p, y\rangle +\langle q, x-y\rangle < \langle p, x\rangle$ for all $(y,q)\in G(T),$}
\]
and therefore 
\[
\hbox{$\langle p-q, x-y\rangle> 0$ for all $(y,q)\in G(T).$}
\] 
But since $T$ is maximal monotone, this means that $p\in Tx$. But then $p\in \bar \partial L_T(x)$ by the first part, leading to $L_T(x,p)= \langle x,p\rangle$, which is a contradiction.

Finally, note that property (\ref{fenchel}) on $L$ yields that  $\bar \partial L(x) \subset \delta L (x)$ and therefore we finally obtain that $Tx=\bar \partial L(x)$.

\begin{proposition} Let  $X$ be a separable reflexive Banach space, and let $L$ be a  convex lower semi-continuous Lagrangian on $X\times X^*$ that satisfies 
 \begin{equation}
\hbox{$L^*(p,x) \geq L(x,p)\geq \langle x, p\rangle$ for every $(x,p)\in X\times X^*$.}
\end{equation}
Then, there exists a selfdual Lagrangian $N$ on $X\times X^*$ such that $\bar \partial L=\bar \partial N$ and 
 \begin{equation}
\hbox{$ L(x,p) \leq N(x,p) \leq L^*(p,x)$ for every $(x,p)\in X\times X^*$.}
 \end{equation}
\end{proposition} 

{\bf Proof:}  Consider the following two Lagrangians $M_0(x,p)=L(x,p)$ and  $N_0(x,p)=L^*(p,x)$ in such a way that  for every $(x,p)\in X\times X^*$, 
\begin{equation}
\hbox{$\langle x, p\rangle \leq M_0(x,p) \leq N_0(x,p)$, }
\end{equation}
\begin{equation}
M^*_0(p,x) = N_0(x,p), 
\end{equation}
and 
\begin{equation}
\hbox{$ \bar \partial M_0(x)= \bar \partial N_0(x)$.}
\end{equation}
We   define by transfinite induction a  family of convex lower semi-continuous Lagrangians $M_\alpha$ and $N_\alpha$ that verify the following properties for every ordinal $\alpha$:
\begin{enumerate}
\item  $M_0\leq M_\alpha \leq M_{\alpha +1} \leq N_{\alpha+1}\leq N_\alpha \leq N_0$. 
\item $M_\alpha^*(p,x)=N_\alpha (x,p)$.
\item $\bar \partial M_\alpha(x)= \bar \partial N_\alpha (x)$.
\end{enumerate}
Starting with $M_0(x,p)$ and  $N_0(x,p)$, we define $(M_\alpha, N_\alpha)_\alpha$  in the  following way:\\
(i) If $\alpha=\beta +1$, then set
 \[
 \hbox{$M_\alpha(x,p)= \inf\left\{\frac{M_\beta (x-y, p-q)}{2} +\frac{N_\beta (x+y, p+q)}{2}; \, (y,q)\in X\times X^*\right\}$}
\]
and 
\[
\hbox{$N_\alpha(x,p)=\frac{1}{2}\big(M_\beta (x,p) +N_\beta (x,p)\big)$.}
\]
 It is then clear that 
\[
M_\beta \leq M_\alpha \leq N_\alpha \leq N_\beta, 
\] 
and that $M_\alpha ^*(p,x)=(\frac{M_\beta}{2}\star \frac{N_\beta}{2})^*(2p,2x)=\big((\frac{M_\beta}{2})^*+(\frac{N_\beta}{2})^*\big)(p,x)=\big(\frac{N_\beta}{2}+ \frac{M_\beta}{2}\big)(x,p)=N_\alpha (x,p)$. \\

(ii) If $\alpha$ is a limit ordinal, then we define:
\[
\hbox{ $M_\alpha(x,p)=\sup\limits_{\beta <\alpha}\ M_\beta(x,p)=\lim\limits_{\beta_n \to\alpha}\uparrow M_{\beta_n}(x,p)$\quad  and \quad $N_\alpha(x,p)=\inf\limits_{\beta < \alpha} N_\beta(x,p)=\lim\limits_{\beta_n \to \alpha}\downarrow N_{\beta_n}(x,p)$.  }
\]
It is again clear that  $M_\beta \leq M_\alpha  \leq  N_\alpha \leq N_\beta$ for every $\beta <\alpha$, and that $M_\alpha^*(p,x)=N_\alpha (x,p)$.

Since $X$ is separable, there exists $\gamma <\Omega$ the first uncountable ordinal such that $M_{\gamma +1}=M_\gamma$, which means that $M_\gamma=N_\gamma$. Set $N(x,p)=M_\gamma (x,p)$ and note that 
\[
N^*(p.x)=M^*_\gamma (p,x)=N_\gamma (x,p)=M_\gamma (x,p)=N(x,p), 
\]
which means that $N$ is a selfdual Lagrangian. Moreover, since for every $\alpha <\Omega$ we have
\[
N_0 (x,p)-\langle x,p\rangle \geq N_\alpha (x,p)-\langle x,p\rangle \geq M_\alpha (x,p)-\langle x,p\rangle \geq M_0 (x,p)-\langle x,p\rangle \geq 0,
\]
it follows that $\bar \partial M_0(x)= \bar \partial N_L (x)=\bar \partial N_0 (x)$.\\

{\bf End of proof of Theorem \ref{ghoussoub}:} Associate to the maximal monotone operator $T$ the sub-selfdual Lagrangian $L_T$ via Lemma \ref{fitz}, that is 
\[
\hbox{$T=\bar \partial L_T$ and $L_T^*(p, x)\geq L_T(x,p)\geq \langle x,p\rangle$.}
\]
Now apply the preceding Proposition to $L_T$ to find a selfdual Lagrangian $N_T$ such that 
 $Tx= \bar \partial N_T(x)$ for any $x\in  {\rm Dom} (T)$.

\end{document}